\documentclass[dvipsnames,a4paper,twoside]{article}
\usepackage{amsthm,amsmath,amssymb,xcolor,soul,xspace}
\usepackage[bookmarks=true,bookmarksopen=true]{hyperref}

\usepackage[ruled,algo2e]{algorithm2e}

\theoremstyle{plain}
\newtheorem{theorem}{Theorem}[section]
\newtheorem{lemma}[theorem]{Lemma}
\newtheorem{corollary}[theorem]{Corollary}

\theoremstyle{definition}
\newtheorem{definition}[theorem]{Definition}

\theoremstyle{remark}
\newtheorem{remark}[theorem]{Remark}

 \numberwithin{equation}{section}

\title{Convergence analysis of the semismooth Newton method for sparse control problems governed by semilinear elliptic equations\thanks{The authors were partially supported by MCIN/ AEI/10.13039/501100011033 under research project PID2020-114837GB-I00.}}

\author{Eduardo Casas\thanks{Departamento de Matem\'{a}tica Aplicada y Ciencias de la Computaci\'{o}n, E.T.S.I. Industriales y de Telecomunicaci\'on, Universidad de Cantabria, 39005 Santander, Spain
(\texttt{eduardo.casas@unican.es})},
\and Mariano Mateos\thanks{Departamento de Matem\'{a}ticas, Campus de Gij\'on, Universidad de Oviedo, 33203, Gij\'on, Spain(\texttt{mmateos@uniovi.es})}
}

\ifpdf
\hypersetup{
  pdftitle={Convergence analysis of the semismooth Newton method for sparse control problems governed by semilinear elliptic equations},
  pdfauthor={E.~Casas, and M.~Mateos}
}
\fi

\newcommand{\dx}{\,\mathrm{d}x}

\newcommand{\Pb}{\mbox{\rm (P)}\xspace}

\newcommand{\uad}{U_{\rm ad}}

\newcommand{\proj}{\operatorname{Proj}}

\begin{document}

\maketitle

\begin{abstract}
We show that a second order sufficient condition for local optimality, along with a strict complementarity condition,
is enough to get the super-linear convergence of the semismooth Newton method for an optimal control problem governed by a semilinear elliptic equation. The objective functional may include a sparsity promoting term and we allow for box control constraints.
\end{abstract}

\begin{quote}
\end{quote}

\begin{quote}
\end{quote}

\section{Introduction}
\label{S1}
Let us consider a domain $\Omega\subset\mathbb R^n$, $n\leq 3$, with a Lipschitz boundary $\Gamma$. We study the following problem:
\[
\Pb\qquad \min_{u \in \uad}  J(u):=F(u)+\gamma j(u)
\]
where
\[F(u) = \int_\Omega L(x,y_u(x)) \dx+ \frac{\kappa}{2}\int_\Omega u(x)^2\dx\text{ and }j(u)= \int_\Omega|u(x)|\dx.\]
Here $L:\Omega\times \mathbb{R}\to\mathbb{R}$ is a given function, $\kappa > 0$, $\gamma \geq 0$ and
\[
\uad =\{u \in L^2(\Omega): \alpha \le u(x) \le \beta\text{ for a.a. } x \in \Omega\},
\]
with $-\infty \le \alpha < \beta \le \infty$. If $\gamma >0$, we will further suppose $\alpha < 0 <\beta$.

Above $y_u$ denotes the state associated to the control $u$ related by the following semilinear elliptic state equation
\begin{equation}\label{E1.1}
\left\{\begin{array}{rcll}
A y_u + f(x,y_u) &=& u&\mbox{ in }\Omega,\\
 y_u& =& 0 &\mbox{ on }\Gamma.
  \end{array}
  \right.
\end{equation}
Assumptions on the data $A$, $f$, $L$ are specified in Section \ref{S2}.

To introduce the main result of the paper and put it in the context of related results in the literature, we briefly describe the semismooth Newton method; precise definitions will be introduced in Section \ref{S3}. Let $\bar u$ be a solution of the equation $\Phi(u)=0$, where $\Phi$ is a semismooth function. Given $u_k$, at every step we select $M_k\in\partial\Phi(u_k)$, the generalized derivative of $\Phi$ at $u_k$, we compute the solution of the linear system $M_k v_k=-\Phi(u_k)$, and set $u_{k+1}=u_k+v_k$. We have that $u_k$ converges superlinearly to  $\bar u$ provided that $u_0$ is close enough to $\bar u$ and the inverses of the operators $M_k$ exist and are uniformly bounded.

In the case of linear equations and convex objective functionals, the uniform boundedness is obtained assuming the  existence of $\nu>0$ such that $F''(\bar u) v^2\geq \nu\|v\|^2_{L^2(\Omega)}$ for all $v\in L^2(\Omega)$; see \cite{Hintermuller-Ito_Kunisch2002,Pieper2015,Stadler2009}.
While this assumption is fully justified in that case,  it is too restrictive if the equation is not linear because it is too far from the second order necessary condition $F''(\bar u)\geq0$ for all $v\in C_{\bar u}$, the cone of critical directions.
As far as we know, the only papers dealing with the convergence of semismooth Newton's method for optimal control problems governed by nonlinear equations are \cite{Amstutz-Laurain}, \cite{Mannel-Rund2021} and \cite{Pieper2015}. In the last two references, the proof of the convergence is done assuming the above mentioned condition, while in the first one a condition implying convexity of the functional is done.

The goal, and the novelty, of our paper is the proof of the superlinear convergence of the semismooth Newton method towards a local solution $\bar u$ of $\Pb$ assuming a strict complementarity  condition, to be properly established in Definition \ref{D2.11}, along with a sufficient second order condition for local optimality. The sufficient second order condition is the usual one enjoying a minimal gap with respect to the necessary one.
In Theorem \ref{T3.4}, we prove that these two hypotheses imply the uniform boundedness of the inverses of the selected generalized derivatives. A strict complementary assumption together with a second order sufficient condition are the usual hypotheses to prove the super-linear convergence of numerical algorithms in finite dimensional constrained optimization problems, cf. \cite{Robinson1972,Robinson1974}; see also \cite[Chapter 17 and 18]{Nocedal-Wright1999} and the references therein.

The plan of the paper is as follows. In Section \ref{S2} we introduce the assumptions on the control problem and carry out the first and second order analysis. The convergence of the semismooth Newton algorithm is proved in Section \ref{S3}. In the last section, we describe some computational details and present two numerical examples.

\section{Assumptions and first and second order analysis of the control problem}
\label{S2}

We make the following assumptions on the data of the control problem.

(A1) Throughout all the paper, $\Omega$ is a bounded open subset of $\mathbb R^n$, $1 \le n\leq 3$. If $n = 2$ or $3$ we assume that its boundary $\Gamma$ is Lipschitz. If $n = 1$, $\Omega$ is a bounded interval and $\Gamma$ is reduced to the two end points of the interval. The operator $A$ is defined in $\Omega$ by the expression
\[
Ay = -\sum_{i, j =1}^n\partial_{x_j}[a_{ij}\partial_{x_i}y] + a_0y
\]
with $a_0, a_{i,j}\in L^\infty(\Omega)$ for $1 \le i, j \le n$, $a_0 \ge 0$, and there exists $\Lambda>0$ such that
\[
\sum_{i,j=1}^n a_{i,j}(x)\xi_i\xi_j\geq \Lambda |\xi|^2\mbox{ for a.e. }x\in\Omega\mbox{ and all }\xi\in\mathbb{R}^n.
\]

(A2) We assume that $f:\Omega\times\mathbb{R}\longrightarrow\mathbb{R}$ is a Carath\'eodory function of class $C^2$ with respect to the second variable satisfying the following conditions for almost all $x \in \Omega$:
\begin{align*}
&\bullet \exists \bar p>\frac{n}{2} \text{ such that } f(\cdot,0) \in L^{\bar p}(\Omega),\\
&\bullet\frac{\partial f}{\partial y}(x,y)\geq 0\quad \forall y\in\mathbb{R},\\
&\bullet\forall M>0,\ \exists C_M >0 \text{ such that } \sum_{j=1}^2\left|\frac{\partial^j f}{\partial y^j}(x,y)\right|\leq C_{f,M} \mbox{ for all }|y|\leq M,\\
&\bullet\forall\varepsilon > 0\mbox{ and }\forall M>0\  \exists \delta >0 \text{ such that }\left|\frac{\partial^2 f}{\partial y^2}(x,y_1) - \frac{\partial^2 f}{\partial y^2}(x,y_2)\right| \le \varepsilon\\
&\qquad\ \ \, \mbox{ for all }|y_1|,|y_2| \le M\mbox{ with }|y_1-y_2|\le\delta.
\end{align*}

(A3) For the cost functional we suppose that $L:\Omega \times\mathbb{R}\longrightarrow\mathbb{R}$ is a Carath\'eodory function of class $C^2$ with respect to the second variable satisfying the following conditions for almost all $x \in \Omega$:
\begin{align*}
&\bullet L(\cdot,0)\in L^1(\Omega) \text{ and } \forall M>0\ \exists \Psi_{L,M}\in L^{\bar p}(\Omega) \text{ and } C_{L,M}>0  \text{ such that }\\
&\qquad\qquad\left|\displaystyle\frac{\partial L}{\partial y}(x,y)\right|\leq \Psi_{L,M}(x)\ \mbox{ and }\
\left|\displaystyle\frac{\partial^2 L}{\partial y^2}(x,y)\right|\leq C_{L,M}
\text{ for all }|y|\leq M,\\
&\bullet\forall\varepsilon > 0\mbox{ and }\forall M>0,\  \exists \delta >0 \text{ such that }\ \left|\displaystyle\frac{\partial^2 L}{\partial y^2}(x,y_1)-
  \frac{\partial^2 L}{\partial y^2}(x,y_2)\right|\le \varepsilon\\
& \qquad\qquad\! \mbox{ for all }|y_1|,|y_2| \le M \mbox{ with }|y_1-y_2|<\delta.
\end{align*}

Let us consider the Banach space $Y = H_0^1(\Omega) \cap C(\bar\Omega)$.
Under the above assumptions, the following properties are well known; see, for instance, \cite[Theorem 1.1.2]{CKT2022}.
\begin{theorem}
For any $u\in L^p(\Omega)$ with $p>n/2$, there exists a unique solution of \eqref{E1.1} $y_u\in Y$. Moreover, there exists a constant $K>0$, that depends on $A$, $\Omega$, $p$ and $\bar p$ such that
\[
\|y_u\|_{H^1_0(\Omega)}+\|y_u\|_{C(\bar\Omega)}\leq K(\|u\|_{L^p(\Omega)}+\|f(\cdot,0)\|_{L^{\bar p}(\Omega)})
\]
holds.
The mapping $S:L^p(\Omega)\longrightarrow Y$ given by $S(u)=y_u$ is of class $C^2$. Furthermore, for all $u, v\in L^p(\Omega)$, $z_v=S'(u)v$ is the unique solution to
\[
  \left\{\begin{array}{rcll}
Az + \displaystyle\frac{\partial f}{\partial y}(x,y_u)z  & = & v&\mbox{ in }\Omega,\\
 z& =& 0 &\mbox{ on }\Gamma,\end{array}\right.
\]
and, given $v_1,v_2\in L^p(\Omega)$, $z_{v_1,v_2}=S''(u)(v_1,v_2)$ is the unique solution to
\[
\left\{\begin{array}{rcll}
Az + \displaystyle\frac{\partial f}{\partial y}(x,y_u)z  & = & -\displaystyle\frac{\partial^2 f}{\partial y^2}(x,y_u)z_{v_1}z_{v_2}&\mbox{ in }\Omega,\\
 w& =& 0 &\mbox{ on }\Gamma,\end{array}\right.
\]
where $z_{v_i} = S'(u)v_i$, $i = 1, 2$.
\label{T2.1}
\end{theorem}
For later reference, it will be useful to define the adjoint state in the following form. We consider the mapping $T:L^\infty(\Omega) \longrightarrow Y$ such that $\varphi = T(y)$ is the unique solution to the adjoint state equation:
\[
\left\{\begin{array}{rcll}
A^*\varphi + \displaystyle\frac{\partial f}{\partial y}(x,y)\varphi  & = & \displaystyle\frac{\partial L}{\partial y}(x,y)&\mbox{ in }\Omega,\\
 \varphi& =& 0 &\mbox{ on }\Gamma.
  \end{array}
  \right.
\]
Setting $G = T \circ S$, we have that the adjoint state related to $u$ is given by $\varphi_u=G(u)$. From Theorem \ref{T2.1} and  the chain rule, it is straightforward to deduce the following two results.
\begin{theorem}\label{T2.2}
  For every $p>n/2$, the mapping $G:L^p(\Omega)\to Y$ is of class $C^1$ and for every $u, v\in L^p(\Omega)$, $\eta_v = G'(u) v$ is the unique solution of
\begin{equation}\label{E2.1}
\left\{\begin{array}{rcll}
A^*\eta_v + \displaystyle\frac{\partial f}{\partial y}(x,y_u)\eta_v   & = & \displaystyle
\left(
\frac{\partial^2 L}{\partial y^2}(x,y_u)
-\varphi_u \frac{\partial^2 f}{\partial y^2}(x,y_u)\right)z_v
&\mbox{ in }\Omega,\\
 \eta_v& =& 0 &\mbox{ on }\Gamma.
  \end{array}
  \right.
\end{equation}
\end{theorem}
\begin{theorem}\label{T2.3}The functional $F:L^2(\Omega)\to\mathbb{R}$ is of class $C^2$. Further, for every $u, v,v_1,v_2\in L^2(\Omega)$ the following identities hold
\begin{align}
F'(u)v = &\int_\Omega(\varphi_u+\kappa u)v\dx,\label{E2.2}\\
F''(u)(v_1,v_2) = &\int_\Omega
\left\{ \left( \frac{\partial^2 L}{\partial y^2}(x,y_u)-\varphi_u \frac{\partial^2 f}{\partial y^2}(x,y_u)\right)z_{v_1} z_{v_2}+\kappa v_1 v_2\right\} \dx\label{E2.3}\\
= & \int_\Omega(\eta_{v_1}+\kappa v_1) v_2\dx =\int_\Omega(\eta_{v_2}+\kappa v_2) v_1\dx,\label{E2.4}
\end{align}
where $\varphi_u = G(u)$,  $z_{v_i}= G'(u)v_i$, and $\eta_{v_i} = G'(u)v_i$ for  $i=1,2$.
\end{theorem}
We will need some results about the adjoint states, which we gather in the next lemma.
\begin{lemma}
Given $R>0$ arbitrary, we denote by $\bar B_R(0)$ the closed $L^2(\Omega)$-ball centered at $0$ with radius $R$. There exist a constant $K_{G'}(R)>0$ such that
\begin{align}
&\|G'(u)v\|_Y\leq K_{G'}(R) \|v\|_{L^2(\Omega)}\quad \forall u \in \bar B_R(0) \text{ and }\forall v \in L^2(\Omega),\label{E2.5}\\
&\|G(u_1)-G(u_2)\|_Y\leq K_{G'}(R)\|u_1-u_2\|_{L^2(\Omega)} \quad \forall u_1, u_2 \in \bar B_R(0).\label{E2.6}
\end{align}
\label{L2.4}
\end{lemma}
\begin{proof}
Let us prove \eqref{E2.5}. From Theorem \ref{T2.1}, we deduce the existence of a constant $M(R)>0$ such that $\|y_u\|_{C(\bar\Omega)}\leq M(R)$ for every $u \in \bar B_R(0)$. Moreover, from the monotonicity of $f$ we deduce the existence of a constant $C_1$ such that
\begin{equation}
\|z_v\|_Y \le C_1\|v\|_{L^2(\Omega)} \quad \forall u \in \bar B_R(0) \text{ and } \forall v \in L^2(\Omega).
\label{E2.7}
\end{equation}
We also obtain with assumption (A3)
\begin{equation}
\|\varphi_u\|_Y \le C_1\Big\|\frac{\partial L}{\partial y}(\cdot,y_u)\Big\|_{L^2(\Omega)} \le C_1\|\Psi_{L,M(R)}\|_{L^2(\Omega)} \le C_R \ \forall u \in \bar B_R(0).
\label{E2.8}
\end{equation}
Once again, from the equation \eqref{E2.1} and using \eqref{E2.7} and \eqref{E2.8} along with the assumptions (A2) and (A3) we get
\[
\|G'(u)v\|_Y = \|\eta_v\|_{Y}\leq K_{G'}(R)\|v\|_{L^2(\Omega)} \quad \forall u \in \bar B_R(0) \text{ and } \forall v \in L^2(\Omega).
\]
Thus, \eqref{E2.5} follows. Estimate \eqref{E2.6} is readily deduced from \eqref{E2.5} and the generalized mean value theorem.
\end{proof}

Let us also remark that $j(u)=\|u\|_{L^1(\Omega)}$ is convex and Lipschitz. For every $u,v\in L^1(\Omega)$, the directional derivative $j'(u;v)$ is given by
\begin{equation}
j'(u;v)= \int_{\Omega_u^+}v\dx-\int_{\Omega_u^-}v\dx+
\int_{\Omega_u^0}|v|\dx,
\label{E2.9}
\end{equation}
where $\Omega_u^+$, $\Omega_u^-$ and $\Omega_u^0$ are the sets of points where $u$ is respectively positive, negative or zero. We denote $J'(u;v) = F'(u)v+\gamma j'(u;v)$ for every $v \in L^2(\Omega)$. The subdifferential of $j$ at $u$ is given by
\begin{equation}\label{E2.10}
\partial j(u) = \left\{\lambda \in L^\infty(\Omega) : \lambda(x)\in\left\{\begin{array}{cc}
                        \{+1\} & \mbox{ if }u(x)>0, \\
                        \{-1\} & \mbox{ if }u(x)<0, \\\relax
                        [-1,1] & \mbox{ if } u(x)=0.
                      \end{array}
                \right.\right\}
\end{equation}
A local solution of \Pb is intended in the $L^2(\Omega)$-sense along this paper. In the following theorem, we summarize necessary and sufficient conditions for local optimality. First, we define the cone of critical directions by
\[
C_{\bar u}=\{v\in L^2(\Omega) : \text{satisfying \eqref{E2.11} and } J'(\bar u)v+\gamma j'(\bar u;v)=0\},
\]
where
\begin{equation}\label{E2.11}
\left\{\begin{array}{l}v(x)\geq 0\mbox{ if }\bar u(x)=\alpha,\\v(x)\leq 0\mbox{ if }\bar u(x)=\beta.\end{array}\right.
\end{equation}

\begin{theorem}\label{T2.5}
  Suppose $\bar u\in\uad$ is a local solution of \Pb. Then,  the following conditions hold
\begin{align}
& J'(\bar u;u-\bar u)\geq 0\ \forall u\in\uad,\label{E2.12}\\
&\exists \bar\lambda\in \partial j(\bar u) \text{ such that } \int_\Omega (\bar\varphi+\kappa\bar u+\gamma\bar\lambda)(u-\bar u)\geq 0\ \forall u\in\uad,\label{E2.13}\\
&F''(\bar u)v^2\geq 0\ \forall v\in C_{\bar u}, \label{E2.14}
\end{align}
where $\bar\varphi =  G(\bar u)$. Conversely, suppose that $(\bar u,\bar\lambda)\in \uad \times \partial j(\bar u)$ satisfies \eqref{E2.13} and
\begin{equation}
\label{E2.15}
F''(\bar u) v^2>0\text{ for all }v\in C_{\bar u}\setminus\{0\}.
\end{equation}
Then, there exist $\nu>0$ and $\delta >0$ such that
\[
J(\bar u) +\frac{\nu}{2}\|u-\bar u\|_{L^2(\Omega)}^2 \leq J(u)\quad \forall u\in \uad \text{ with } \|u-\bar u\|_{L^2(\Omega)} \le \delta.
\]
\end{theorem}

 The reader is referred to \cite[Theorems 3.1, 3.7 and 3.9]{CHW20212-SIOPT} for its proof.
Notice that the gap between the sufficient condition \eqref{E2.15} and the necessary condition \eqref{E2.14} is the minimal one, the same as in finite dimensional optimization.

We quote the following result, whose proof can be found in \cite[Corollary 3.2]{CHW20212-SIOPT}.
\begin{corollary}\label{C2.6}
Let $(\bar u,\bar\varphi,\bar\lambda) \in\uad \times Y \times \partial j(\bar u)$ satisfy \eqref{E2.13} with $\bar\varphi = G(\bar u)$. Then, the following relation holds:
\begin{equation}
\bar u(x) = \proj_{[\alpha,\beta]}\left(-\frac{1}{\kappa}\big(\bar\varphi(x) + \gamma\bar\lambda(x)\big)\right). \label{E2.16}
\end{equation}
Moreover, if $\gamma > 0$ the following properties are fulfilled:
 \begin{align}
     \bar u(x)  =& 0\iff |\bar\varphi(x)|\leq\gamma,\label{E2.17}\\
    \bar\lambda(x) = & \proj_{[-1,+1]}\left(-\frac{1}{\gamma}\varphi(x)\right).\label{E2.18}
 \end{align}

\end{corollary}

\begin{remark}\label{R2.7}
Notice that if $\gamma=0$, the role of $\bar\lambda$ in \eqref{E2.13} and \eqref{E2.16} is irrelevant. Since, nevertheless, the notation is consistent, we leave it there in order to make an exposition as unified as possible of both cases, $\gamma=0$ and $\gamma>0$.
\end{remark}

\begin{remark}\label{R2.8}
As immediate consequence of \eqref{E2.16} we obtain
\begin{align}
&\text{If } \bar u(x) = \alpha\ \text{ then } \ \bar\varphi(x)+\kappa\bar u(x)+\gamma\bar\lambda(x) \geq 0.  \label{E2.19}\\
&\text{If } \bar u(x) = \beta\ \text{ then } \ \bar\varphi(x)+\kappa\bar u(x)+\gamma\bar\lambda(x) \leq 0.  \label{E2.20}\\
&\text{If } \alpha < \bar u(x) < \beta \ \text{ then } \ \bar\varphi(x)+\kappa\bar u(x)+\gamma\bar\lambda(x) = 0.\label{E2.21}\\
&\text{If } \bar\varphi(x)+\kappa\bar u(x)+\gamma\bar\lambda(x) > 0  \ \text{ then } \ \bar u(x) = \alpha.\label{E2.22}\\
&\text{If } \bar\varphi(x)+\kappa\bar u(x)+\gamma\bar\lambda(x) < 0  \ \text{ then } \ \bar u(x) = \beta.\label{E2.23}
\end{align}
Using \eqref{E2.2}, \eqref{E2.9}, and \eqref{E2.10} we infer
\begin{equation}
J'(\bar u;v) = \int_{\Omega_{\bar u}^+\cup\Omega_{\bar u}^-} [\bar\varphi(x)+\kappa\bar u(x)+\gamma\bar\lambda(x)]v(x)\dx + \int_{\Omega^0_{\bar u}}[\bar\varphi(x)v(x) + \gamma|v(x)|]\dx.
\label{E2.24}
\end{equation}
\end{remark}

The next lemma establishes an important property of the elements of the critical cone.

\begin{lemma}\label{L2.9}
Let $(\bar u,\bar\varphi,\bar\lambda)$ be as in Corollary \ref{C2.6}. Then, the following property holds for $v \in L^2(\Omega)$ and for almost all $x \in \Omega$
\begin{equation}
[\bar\varphi(x) + \kappa\bar u(x) + \gamma\bar\lambda(x)]v(x)\left\{\begin{array}{ll}\ge 0&\text{if $v$ satisfies \eqref{E2.11}},\\= 0&\text{if } v \in C_{\bar u}.\end{array}\right.
\label{E2.25}
\end{equation}
\end{lemma}

\begin{proof}
The inequality of \eqref{E2.25} is a straightforward consequence of \eqref{E2.19}--\eqref{E2.21} and \eqref{E2.11}. To prove the equality of \eqref{E2.25} we recall that $j'(\bar u;v) \ge \int_\Omega\lambda v \dx$ for all $\lambda \in \partial j(\bar u)$. Then, we get
\[
0 = J'(\bar u;v) \ge \int_\Omega[\bar\varphi(x) + \kappa\bar u(x) + \gamma\bar\lambda(x)]v(x)\dx\quad \forall v \in C_{\bar u}.
\]
Since the integrand is non-negative for almost all $x \in \Omega$, the above inequality yields $[\bar\varphi(x) + \kappa\bar u(x) + \gamma\bar\lambda(x)]v(x) = 0$ for almost all $x \in \Omega$.
\end{proof}

\begin{lemma}
Let $(\bar u,\bar\varphi,\bar\lambda)$ be as in Corollary \ref{C2.6}. Then, $C_{\bar u}$ is the set of elements $v \in L^2(\Omega)$ satisfying the following conditions:
\begin{align}
&v(x) = 0 \text{ if } \bar\varphi(x) + \kappa\bar u(x) + \gamma\bar\lambda(x) \neq 0 \text{ or } |\bar\varphi(x)| < \gamma, \label{E2.26}\\
&\left\{\begin{array}{l}v(x)\geq 0\mbox{ if }\bar u(x)=\alpha \text{ or } \bar\varphi(x) = - \gamma,\\v(x)\leq 0\mbox{ if }\bar u(x)=\beta \text{ or } \bar\varphi(x) = + \gamma,\end{array}\right.
\label{E2.27}
\end{align}
where the terms involving $\gamma$ should be removed if $\gamma = 0$.
\label{L2.8}
\end{lemma}

\begin{proof}
From \eqref{E2.17} we infer that $\bar\varphi(x)v(x) + \gamma|v(x)| \ge 0$ for almost all $x \in \Omega^0_{\bar u}$. Using this, the inequality in \eqref{E2.25}, and \eqref{E2.24} we deduce that $J'(\bar u;v) = 0$ if and only if $[\bar\varphi(x) + \kappa\bar u(x) + \gamma\bar\lambda(x)]v(x) = 0$ a.e.~in $\Omega_{\bar u}^+\cup\Omega_{\bar u}^-$ and $\bar\varphi(x)v(x) + \gamma|v(x)| = 0$ a.e.~in $\Omega^0_{\bar u}$. The last equality holds if and only if ($v(x) = 0$ if $|\bar\varphi(x)| < \gamma$), ($v(x) \ge 0$ if $\bar\varphi(x) = -\gamma$), and ($v(x) \le 0$ if $\bar\varphi(x) = +\gamma$). These equivalences  prove the  characterization of $C_{\bar u}$ given in the statement of the lemma.
\end{proof}

Now, we define the following closed vector subspace of $L^2(\Omega)$
\[
T_{\bar u} = \{v\in L^2(\Omega):\ v(x)=0\mbox{ if }|\bar\varphi(x)+\kappa\bar u(x)+\gamma\bar\lambda(x)|>0\ \text{ or }\ |\bar\varphi(x)| < \gamma\}
\]
and the set
\[
\Sigma_{\bar u} = \{x\in\Omega:\ (\bar u(x)\in\{\alpha,\beta\}\mbox{ and }\bar\varphi(x)+\kappa\bar u(x)+\gamma\bar\lambda(x)=0) \text{ or } |\bar\varphi(x)|=\gamma\}.
\]
Once again, the terms involving $\gamma$ should be removed in the case $\gamma = 0$.

\begin{definition}\label{D2.11}
We say that the strict complementary condition is satisfied at $\bar u$ if $|\Sigma_{\bar u}| = 0$, where $| \cdot |$ stands for the Lebesgue measure.
\end{definition}

\begin{lemma}\label{L2.12}
Let $(\bar u,\bar\varphi,\bar\lambda)$ be as in Corollary \ref{C2.6} and assume that the strict complementary condition holds at $\bar u$, then $C_{\bar u}=T_{\bar u}$.
\end{lemma}

This lemma is an immediate consequence of Lemma \ref{L2.8} and the fact that ${|\Sigma_{\bar u}| = 0}$.

Given $\tau >0$, where $\tau < \gamma$ if $\gamma > 0$, we define the extended subspace
\[
T^\tau_{\bar u} = \{v\in L^2(\Omega): v(x)=0\mbox{ if }|\bar\varphi(x)+\kappa\bar u(x)+\gamma\bar\lambda(x)|>\tau\ \text{ or }\ |\bar\varphi(x)| < \gamma - \tau\}.
\]

\begin{theorem}\label{T2.13}
Let $(\bar u,\bar\varphi,\bar\lambda) \in\uad \times Y \times \partial j(\bar u)$ satisfy \eqref{E2.13} with $\bar\varphi = G(\bar u)$ and assume that the strict complementary condition $|\Sigma_{\bar u}| = 0$ and the second order sufficient condition \eqref{E2.15} hold at $\bar u$.
Then, there exist $\nu >0$ and $\tau > 0$, with $\tau<\gamma$ if $\gamma >0$, such that
\begin{equation}
\label{E2.28}
F''(\bar u)v^2\geq \nu\|v\|^2_{L^2(\Omega)}\text{ for all } v\in  T^\tau_{\bar u}.
\end{equation}
\end{theorem}
\begin{proof}
 We will proceed by contradiction: suppose \eqref{E2.28} is false. Then, there exists a sequence $\{v_k\}_{k=1}^\infty\subset L^2(\Omega)$ such that $v_k\in T_{\bar u}^{1/k}$ and $F''(\bar u)v_k^2<\frac{1}{k}\|v_k\|^2_{L^2(\Omega)}$. Of course, we can assume that $\|v_k\|_{L^2(\Omega)}=1$, otherwise it is enough to divide $v_k$ by its $L^2(\Omega)$-norm to have
\begin{equation}
v_k\in T_{\bar u}^{1/k},\ \|v_k\|_{L^2(\Omega)}=1, \text{ and } F''(\bar u)v_k^2<\frac{1}{k}.
\label{E2.29}
\end{equation}
Then, for a subsequence, denoted in the same way, there exists $v\in L^2(\Omega)$ such that $v_k\rightharpoonup v$ weakly in $L^2(\Omega)$. We observe that $v \in T_{\bar u}$. Indeed, for every $\varepsilon > 0$ we set
\[
\Theta^\varepsilon = \{x\in\Omega : v(x) \neq 0 \ \text{ and }\ |\bar\varphi(x)+\kappa\bar u(x)+\gamma\bar\lambda(x)|>\varepsilon \text{ or }\ |\bar\varphi(x)| < \gamma - \varepsilon\}.
\]
Due to $v_k$ vanishes in $\Theta^\varepsilon$ for every $k > \frac{1}{\varepsilon}$, its weak limit $v$ vanishes as well in $\Theta^\varepsilon$. As $\varepsilon > 0$ is arbitrary, we conclude that $v \in T_{\bar u}$. On the other hand, since the quadratic form $F''(\bar u):L^2(\Omega)\to \mathbb{R}$ is weakly lower semicontinuous, using \eqref{E2.29} we have that
\[
F''(\bar u)v^2\leq \liminf_{k\to\infty} F''(\bar u) v_k^2 = 0.
\]
The strict complementarity condition  implies that $C_{\bar u} = T_{\bar u}$. Therefore, as a consequence of \eqref{E2.15}, we deduce that $v=0$. Moreover, the weak convergence $v_k\rightharpoonup 0$ in $L^2(\Omega)$ implies the strong convergence $z_{v_k}\to 0$ in $C(\bar\Omega)$. Using that $\|v_k\|_{L^2(\Omega)}=1$ we obtain
\[
\lim_{k\to\infty}F''(\bar u)v_k^2 = \lim_{k\to\infty} \int_\Omega\left( \frac{\partial^2 L}{\partial y^2}(x,\bar y)-\bar \varphi \frac{\partial^2 f}{\partial y^2}(x,\bar y)\right)z_{v_k}^2 \dx+\kappa = \kappa,
\]
which contradicts the fact that $\kappa>0$.
\end{proof}

\begin{remark}\label{R2.14}
 In \cite[eq. (3.19)]{Troltz1999}, the author makes an assumption similar to \eqref{E2.28} to prove quadratic convergence for a sequential quadratic programming algorithm. However, \eqref{E2.28} sounds quite strong as an assumption because it seems to be very far from the second order necessary condition. Theorem \ref{T2.13} shows that it is satisfied whenever the no gap second order sufficient condition plus the strict complementarity condition hold.
\end{remark}

\section{Semismooth Newton method}
\label{S3}

Next we use \eqref{E2.16} and \eqref{E2.18} to define an equation $\Phi(u)=0$ satisfied by any local solution of \Pb, where $\Phi$ is semismooth. We define semismoothness following \cite[Definition 3.1]{Ulbrich2011}. A slightly different approach using the concept of slant differentiability can be found in \cite{Hintermuller-Ito_Kunisch2002}.
\begin{definition}\label{D3.1}
   Given two Banach spaces $X$ and $Y$, an open subset $V$ of $X$, a continuous function ${\Phi:V\to Y}$, and a set-valued mapping ${\partial\Phi:V\rightrightarrows \mathcal{L}(X,Y)}$ such that $\partial\Phi(u)\neq\emptyset$ for every $u\in V$, we say that $\Phi$ is $\partial\Phi$-semismooth at $\bar u\in V$ if
  \begin{equation}\label{E3.1}\lim_{v\to 0}\sup_{M\in\partial \Phi(\bar u+v)}\frac{\|\Phi(\bar u+v)-\Phi(\bar u)-Mv\|_{Y}}{\|v\|_{X}} =0.
  \end{equation}

  The multifunction $\partial\Phi$ is called the generalized derivative of $\Phi$.
\end{definition}
The semismooth Newton method spans a sequence according to the Algorithm \ref{Alg1}.

\LinesNumbered
\begin{algorithm2e}[h!]
\caption{Semismooth Newton method.}\label{Alg1}
\DontPrintSemicolon
Initialize Choose $u_0\in V$. Set $j=0$.\;
\Repeat{convergence}{
Choose $M_j\in\partial\Phi(u_j)$ and solve $M_jv_j=-\Phi(u_j)$\label{line3}.\;
Set $u_{j+1}=u_j+v_j$ and $j=j+1$.\;
}
\end{algorithm2e}

The proof of the following convergence theorem can be found in \cite[Theorem 3.13]{Ulbrich2011}. See also \cite[Theorem 1.1]{Hintermuller-Ito_Kunisch2002}.
\begin{theorem}\label{T3.2}
  Suppose that $\Phi:V\to Y$ is  $\partial\Phi$-semismooth at $\bar u\in V$ solution of $\Phi(u)=0$ locally unique. Suppose, furthermore, that the following regularity condition is satisfied: for every $j$, the operator $M_j\in \partial \Phi(u_j)$ is invertible and there exists $C_\Phi>0$ such that
  \begin{equation}\label{E3.2}\|M_{j}^{-1}\|_{\mathcal{L}(Y,X)}\leq C_\Phi\ \forall j\geq 0.\end{equation}
  Then, there exists $\delta >0$ such that for all $u_0\in V$ with $\|u_0-\bar u\|_{X}<\delta$ the sequence $\{u_j\}_{j\geq 0}$ spanned by the semismooth Newton method converges superlinearly to $\bar u$.
\end{theorem}

Taking into account \eqref{E2.16} and \eqref{E2.18}, we< define $\psi:\mathbb{R}\to\mathbb{R}$ as
\begin{align*}
\psi(t) = &
\proj_{[\alpha,\beta]}\left\{
-\frac{1}{\kappa}\left[
t +  \proj_{[-\gamma,+\gamma]}\left(
-t
\right)
\right]
\right\}
\end{align*}
and the superposition operator $\Psi_G:L^2(\Omega)\to L^2(\Omega)$ by $\Psi_G(u)(x) = \psi(G(u)(x))$. We recall that $G(u) = \varphi_u$. We consider the mapping $\Phi:L^2(\Omega)\to L^2(\Omega)$ given by $\Phi(u) = u-\Psi_G(u)$. Corollary \ref{C2.6} implies that $\bar u$ satisfies the equation $\Phi(\bar u)=0$. Next, we study the semismoothness properties of $\Phi$.
Hereafter, $\partial^{\mathrm{CL}}\psi$ denotes the generalized derivative of $\psi$ in the sense of Clarke \cite[Definition 3.10]{ClarkeBook2013}. We observe that $\psi$ is a Lipschitz function.

\begin{lemma}\label{L3.3}The function $\Phi:L^2(\Omega)\to L^2(\Omega)$ is $\partial\Phi$-semismooth at every $u\in L^2(\Omega)$ for the set valued mapping
\[\partial\Phi(u)=\{M=I-N:\ N\in \partial\Psi_G(u)\},\]
where
\begin{align*}
\partial&\Psi_G(u) =
\{N\in \ \mathcal{L}(L^2(\Omega),L^2(\Omega)):\text{there exists a Lebesgue measurable function }h\\
&\text{ such that }
h(x)\in \partial^{\mathrm{CL}}\psi(G(u)(x))\text{ a.e. in }\Omega\text{ and }
Nv = h\cdot G'(u)v\ \forall v\in L^2(\Omega)\}.
\end{align*}
\end{lemma}

\begin{proof}Clearly, $\Phi$ is continuous for being the composition of continuous functions.

It is straightforward to check that the generalized derivative of $\psi$ in the sense of Clarke, see \cite[Theorem 10.27]{ClarkeBook2013}, is given by the expression
\begin{align*}
\partial^{\mathrm{CL}}\psi(t) = & \left\{\begin{array}{cl}
\{0\}                      & \mbox{ if } t\in (-\infty,-\gamma-\kappa\beta)\cup (-\gamma,\gamma)\cup (\gamma-\kappa\alpha,+\infty), \\ \\
\left\{-\frac{1}{\kappa}\right\} & \mbox{ if }t\in (-\gamma -\kappa\beta,-\gamma)\cup(\gamma,\gamma-\kappa\alpha), \\ \\{}
[-\frac{1}{\kappa},0]         & \mbox{ if } t\in \{-\gamma -\kappa\beta,-\gamma,\gamma, \gamma-\kappa\alpha\}.
\end{array}
\right.& & \text{ if }\gamma > 0,\\ \\
\partial^{\mathrm{CL}}\psi(t) = & \left\{\begin{array}{cl}
\{0\}                      & \mbox{ if } t\in (-\infty,-\kappa\beta)\cup  (-\kappa\alpha,+\infty), \\ \\
\left\{-\frac{1}{\kappa}\right\} & \mbox{ if }t\in (-\kappa\beta, -\kappa\alpha), \\ \\{}
[-\frac{1}{\kappa},0]         & \mbox{ if } t\in \{ -\kappa\beta,-\kappa\alpha\}.
\end{array}
\right.& & \text{ if }\gamma = 0.
\end{align*}
Since $\psi$ is piecewise $C^2$, it is 1-order semismooth; see \cite[Proposition 2.26]{Ulbrich2011}.
Thanks to the Lipschitz continuity of $G$, see \eqref{E2.6}, we deduce straightforward from \cite[Theorem 3.49]{Ulbrich2011} that $\Psi_G$ is $\partial\Psi_G$-semismooth. Hence, defining $\partial\Phi(u)=\{M=I-N:\ N\in \partial\Psi_G(u)\}$, we readily obtain that $\Phi$ is $\partial\Phi$-semismooth.
\end{proof}

To perform the step in line \ref{line3} of Algorithm \ref{Alg1}, we have to choose some element in $\partial \Phi(u)$. In order to do this selection and obtain a family of uniformly invertible operators, we define
\begin{align*}
g(t) = & \left\{\begin{array}{cl}
0                     & \mbox{ if } t\in (-\infty,-\gamma-\kappa\beta]\cup [-\gamma,\gamma]\cup [\gamma-\kappa\alpha,+\infty), \\ \\
-\dfrac{1}{\kappa}     & \mbox{ if }t\in (-\gamma -\kappa\beta,-\gamma)\cup(\gamma,\gamma-\kappa\alpha),
\end{array}
\right. &&\text{ if }\gamma > 0,\\ \\
g(t) = & \left\{\begin{array}{cl}
0                     & \mbox{ if } t\in (-\infty,-\kappa\beta] \cup [-\kappa\alpha,+\infty), \\ \\
-\dfrac{1}{\kappa}     & \mbox{ if }t\in (-\kappa\beta,-\kappa\alpha),
\end{array}
\right. &&\text{ if }\gamma = 0.
\end{align*}
Notice that $g(t)=0$ when $t\in \{-\gamma -\kappa\beta,-\gamma,\gamma, \gamma-\kappa\alpha\}$ if $\gamma > 0$ or when  $t\in \{ -\kappa\beta, -\kappa\alpha\}$ if $\gamma=0$ and hence $g(t)\in\partial^{\mathrm{CL}}\psi(t)$.
For a given control $u\in L^2(\Omega)$, we select $M_u\in\partial\Phi(u)$ defined as $M_u v = v - h_u\cdot G'(u)v$, where $h_u(x) = g(\varphi_u(x))$.

\begin{theorem}\label{T3.4}
Let $(\bar u,\bar\varphi,\bar\lambda) \in\uad \times Y \times \partial j(\bar u)$ satisfy \eqref{E2.13} with $\bar\varphi = G(\bar u)$ and assume that the strict complementary condition $|\Sigma_{\bar u}| = 0$ and the second order sufficient condition \eqref{E2.15} hold at $\bar u$.  Then, there exist $\delta >0$ and $C>0$ such that for all $u\in B_\delta(\bar u)$ and all $w\in  L^2(\Omega)$, the equation $M_uv=w$ has a unique solution $v\in L^2 (\Omega)$ and the inequality $\|v\|_{L^2(\Omega)}\leq C \|w\|_{L^2(\Omega)}$ holds.
\end{theorem}
\begin{proof}
We define the active and inactive sets for $u$. For $\gamma >0$ we define
\begin{align*}
  \mathbb{A}^\beta  =& \{x\in\Omega:\  \varphi_u(x) \leq -\gamma-\kappa\beta\}, \\
  \mathbb{J}^+      =& \{x\in\Omega:\ -\gamma-\kappa\beta < \varphi_u(x) < - \gamma\},\ \\
  \mathbb{A}^0      =& \{x\in\Omega:\ |\varphi_u(x)|\leq \gamma\}, \\
  \mathbb{J}^-      =& \{x\in\Omega:\ \gamma < \varphi_u(x)  <  \gamma - \kappa\alpha\}, \\
  \mathbb{A}^\alpha =& \{x\in\Omega:\ \gamma - \kappa\alpha \leq \varphi_u(x)\}.
\end{align*}
Notice that all the five sets are disjoint and their union is $\Omega$. We set $\mathbb{A}=\mathbb{A}^\alpha\cup \mathbb{A}^\beta\cup\mathbb{A}^0$ and $\mathbb{J}=\mathbb{J}^-\cup\mathbb{J}^+$.
In the case $\gamma =0$ we define  $\mathbb{A}=\mathbb{A}^\alpha\cup \mathbb{A}^\beta$ and
\begin{align*}
  \mathbb{J}      =& \{x\in\Omega:\ -\kappa\beta < \varphi_u(x) <  - \kappa\alpha\}.
\end{align*}

Using the notation $\eta_v=G'(u)v$,  we have that
\begin{equation*}
M_u v =
\left\{
\begin{array}{ll}
  v &\mbox{ in }  \mathbb{A},\\
   v+\dfrac{1}{\kappa}\eta_v  &\mbox{ in }  \mathbb{J},
\end{array}
\right.
\end{equation*}
and the equation $M_uv=w$ is equivalent to the system
\begin{equation}\label{E3.3}
\left\{
\begin{array}{rcll}
  v &= & w&\mbox{ in }  \mathbb{A},\\
   v+\dfrac{1}{\kappa}\eta_v  &= & w&\mbox{ in }  \mathbb{J}.
\end{array}
\right.
\end{equation}

We write $v = \chi_{{}_\mathbb{J}} v + \chi_{{}_\mathbb{A}} v$. The first equation determines $v$ in the active set $\mathbb{A}$ and we write the second equation as
\begin{equation}\label{E3.4}
  \chi_{{}_\mathbb{J}} v+\frac{1}{\kappa}\eta_{\chi_{{}_\mathbb{J}}v}  =  w- \frac{1}{\kappa}\eta_{\chi_{{}_\mathbb{A}}w} \mbox{ in }  \mathbb{J}.\end{equation}From \eqref{E2.4} we get that this equation is the optimality condition of the unconstrained quadratic optimization problem
\begin{equation}\label{E3.5}
\min_{v\in L^2(\mathbb{J})} H(v):=\frac{1}{2} F''(u)(\chi_{{}_\mathbb{J}}v)^2- \int_{\mathbb{J}} (\kappa w-\eta_{\chi_{{}_\mathbb{A}}w}) v\dx.
\end{equation}
Therefore, if we prove that $H$ has a unique local minimizer in $L^2(\mathbb{J})$, the existence and  uniqueness of a solution of \eqref{E3.4} follows. Using the continuity of the functional $u\to F''(u)$, we deduce the existence of $\delta_0>0$ such that, if $\|u-\bar u\|_{L^2(\Omega)} < \delta_0$, then $|(F''(u)-F''(\bar u))v^2|<\nu/2\|v\|_{L^2(\Omega)}^2$. From Theorem \ref{T2.13}, we deduce the existence of $\tau>0$, with $\tau<\gamma$ if $\gamma >0$, such that
\begin{equation}\label{E3.6}
F''(u)v^2\geq \frac{\nu}{2}\|v\|_{L^2(\Omega)}^2\ \forall v\in T_{\bar u}^\tau\text{ if }\|u-\bar u\|_{L^2(\Omega)} < \delta_0.
\end{equation}
Therefore, \eqref{E3.5} has a unique local minimizer, that is also global, if $L^2(\mathbb{J})\subset T^\tau_{\bar u}$. This embedding follows from the inclusion
\begin{equation}\label{E3.7}
\mathbb{J}\subset\{x\in\Omega:\ |\bar\varphi(x)+\kappa \bar u(x) +\gamma\bar\lambda(x)|\leq\tau\text{ and } |\bar\varphi(x)|\geq\gamma-\tau\},
\end{equation}
or equivalently
\[
\{x\in\Omega:\ |\bar\varphi(x)+\kappa \bar u(x) +\gamma\bar\lambda(x)| > \tau\text{ or } |\bar\varphi(x)| < \gamma-\tau\}\subset \mathbb{A}.
\]
Let us check this inclusion.
Taking $\delta=\min\{\delta_0,1,\dfrac{\tau}{K_{G'}(\bar R)}\}$ with $\bar R = \|\bar u\|_{L^2(\Omega)}+1$, we deduce from \eqref{E2.6} that $\|\varphi_u-\bar\varphi\|_{C(\bar\Omega)} < \tau$ if $\|u-\bar u\|_{L^2(\Omega)} < \delta$.

\underline{Case 1.}
Suppose $\bar\varphi(x)+\kappa \bar u(x) +\gamma\bar\lambda(x) > \tau$.
From \eqref{E2.22}, we have that $\bar u(x) =\alpha$.
If $\gamma > 0$, we also deduce from \eqref{E2.17} and \eqref{E2.18} that $\bar\lambda(x) = -1$. We can write that
 $\bar\varphi(x) >\tau-\kappa \bar u(x) -\gamma\bar\lambda(x) = \tau-\kappa \alpha +\gamma$. Since $\varphi_u(x)-\bar\varphi(x) > -\tau$, we have that
\[
\varphi_u(x) = \varphi_u(x) -\bar\varphi(x)+\bar\varphi(x) > -\tau + \tau-\kappa \alpha +\gamma = -\kappa \alpha +\gamma
\]
and, hence, $x\in \mathbb{A}^\alpha\subset\mathbb{A}$.

\underline{Case 2.}
Suppose $\bar\varphi(x)+\kappa \bar u(x) +\gamma\bar\lambda(x) < -\tau$.
From \eqref{E2.23}, we have that $\bar u(x) =\beta$.
If $\gamma > 0$, we also deduce from \eqref{E2.17} and \eqref{E2.18} that $\bar\lambda(x) = 1$.We can write that
 $\bar\varphi(x) <-\tau-\kappa\bar u(x) -\gamma\lambda(x) = -\tau-\kappa \beta -\gamma$. Since $\varphi_u(x)-\bar\varphi(x) < \tau$, we have that
\[
\varphi_u(x) = \varphi_u(x) -\bar\varphi(x)+\bar\varphi(x) < \tau - \tau-\kappa \beta -\gamma = -\kappa \beta -\gamma
\]
and, consequently, $x\in \mathbb{A}^\beta\subset\mathbb{A}$.

For $\gamma = 0$, cases 1 and 2 imply \eqref{E3.7}.

\underline{Case 3.} Suppose $\gamma > 0$ and $|\bar\varphi(x)|<\gamma-\tau$. Then $|\varphi_u(x)|\leq |\varphi_u(x)-\bar\varphi(x)|+|\bar\varphi(x)| <\tau+\gamma-\tau = \gamma$, which yields $x\in \mathbb{A}^0\subset\mathbb{A}$.

Therefore \eqref{E3.6} and \eqref{E3.7} hold and consequently the system \eqref{E3.3} has a unique solution $v$ and $\chi_{{}_\mathbb{J}}v\in T_{\bar u}^\tau$. It remains to get an estimate for $v$ in terms of $w$ with a constant independent of $u \in B_\delta(\bar u)$.
Using  \eqref{E3.6}, \eqref{E2.4}, and \eqref{E3.4} we get
\begin{align*}
\frac{\nu}{2}\|\chi_{{}_\mathbb{J}}v\|_{L^2(\Omega)}^2\leq &\, F''(u)(\chi_{{}_\mathbb{J}}v)^2 = \int_\Omega (\eta_{\chi_{{}_\mathbb{J}}v}+\kappa \chi_{{}_\mathbb{J}}v)\chi_{{}_\mathbb{J}}v\mathrm{d}x\\
= & \kappa  \int_\Omega \left(w-\frac{1}{\kappa} \eta_{\chi_{{}_\mathbb{A}}w} \right)\chi_{{}_\mathbb{J}}v\mathrm{d}x.
\end{align*}
From the first equation in \eqref{E3.3}, we have that
\[\|\chi_{{}_\mathbb{A}}v\|_{L^2(\Omega)}^2 = \int_\Omega w \chi_{{}_\mathbb{A}}v\mathrm{d}x.\]
Multiplying this equality by $\kappa$ and adding it to the previous inequality, we obtain with the Cauchy-Schwarz inequality and estimate \eqref{E2.5}
\begin{align*}
\min\{\kappa,\frac{\nu}{2}\}\|v\|_{L^2(\Omega)}^2 \leq  & \kappa \int_\Omega  w v\mathrm{d}x - \int_\Omega \eta_{\chi_{{}_\mathbb{A}}w} \chi_{{}_\mathbb{J}}v\mathrm{d}x\\
\leq & \left(\kappa \|w\|_{L^2(\Omega)}  + \|\eta_{\chi_{{}_\mathbb{A}}w}\|_{L^2(\Omega)}\right) \|v\|_{L^2(\Omega)}\\
\leq & \left(\kappa+ K_{G'}(\bar R)\right) \|w\|_{L^2(\Omega)} \|v\|_{L^2(\Omega)}.
\end{align*}
This yields $\|v\|_{L^2(\Omega)} \le C\|w\|_{L^2(\Omega)}$ with $C = \frac{\kappa+ K_{G'}(\bar R)}{\min\{\kappa,\frac{\nu}{2}\}}$.
\end{proof}

The following result is an immediate consequence of Theorem \ref{T3.2}, Lemma \ref{L3.3} and Theorem \ref{T3.4}.
\begin{corollary}\label{C3.5}
Let $(\bar u,\bar\varphi,\bar\lambda) \in\uad \times Y \times \partial j(\bar u)$ satisfy \eqref{E2.13} with $\bar\varphi = G(\bar u)$ and assume that the strict complementary condition $|\Sigma_{\bar u}| = 0$ and the second order sufficient condition \eqref{E2.15} hold at $\bar u$. Then, there exists $\delta >0$ such that for all $u_0\in B_\delta(\bar u)$, the sequence spanned by Algorithm \ref{Alg2} converges superlinearly to $\bar u$.
\end{corollary}

Semismooth Newton's method for problem \Pb is detailed in Algorithm \ref{Alg2}.

\begin{remark}\label{R3.6}Since $C_{\bar u}\subset T_{\bar u}^\tau$, then \eqref{E3.6} implies the sufficient second order optimality condition \eqref{E2.15}.
The proof of Theorem \ref{T3.4} uses \eqref{E3.6}, but the strict complementarity condition $|\Sigma_{\bar u}|=0$ is not used.
Consequently, the statement of Theorem \ref{T3.4}, and also that of Corollary \ref{C3.5}, can be rewritten replacing \eqref{E2.15} and the strict complementarity condition by \eqref{E3.6}. Obviously, this is is a weaker assumption, but it is a less natural assumption: if the strict complementarity condition is not satisfied, the gap between \eqref{E3.6} and the second order necessary condition is too large.
\end{remark}

\begin{algorithm2e}[h!]
\caption{Semismooth Newton method to solve \Pb.}\label{Alg2}
\DontPrintSemicolon
Initialize Choose $u_0\in L^2(\Omega)$. Set $j=0$.\;
\Repeat{convergence}{
Compute $y_j = S(u_j)$ solving the nonlinear equation\label{line3alg2}
\[A y_j + f(x,y_j) = u_j\text{ in }\Omega,\ y_j=0\text{ in }\Gamma\]
\;
Compute $\varphi_j = G(u_j)$ solving the linear equation
\[A^*\varphi_j +\frac{\partial f}{\partial y}(x,y_j)\varphi_j = \frac{\partial L}{\partial y}(x,y_j)\text{ in }\Omega,\ \varphi_j=0\text{ in }\Gamma\]
\;
Compute $\mathbb{A}^\beta_j$, $\mathbb{A}^0_j$, $\mathbb{A}^\alpha_j$, $\mathbb{A}_j$, and  $\mathbb{J}^+_j$,  $\mathbb{J}^-_j$, $\mathbb{J}_j$ using $\varphi_j$. \;
Compute
\[w_j(x) = -\Phi(u_j)(x) = \left\{\begin{array}{lc}
                              -u_j(x) +\beta & \text{ if }x \in \mathbb{A}^\beta_j \\
                              -u_j(x) - \frac{1}{\kappa}(\varphi_k(x)+\gamma) & \text{ if }x \in \mathbb{J}^+_j \\
                              -u_j(x)  & \text{ if }x \in \mathbb{A}^0_j \\
                              -u_j(x) - \frac{1}{\kappa}(\varphi_k(x)-\gamma) & \text{ if }x \in \mathbb{J}^-_j \\
                              -u_j(x) +\alpha & \text{ if }x \in \mathbb{A}^\alpha_j
                            \end{array}\right.
\]
\;
Compute $\eta_j = \eta_{\chi_{{}_{\mathbb{A}_j}}w_{j}}$ solving the linear equations\label{line8}
\begin{align*}
A z_j +\frac{\partial f}{\partial y}(x,y_j)z_j =& \chi_{{}_{\mathbb{A}_j}}w_{j}\text{ in }\Omega,\ z_j = 0\text{ on }\Gamma\\
A^*\eta_j +\frac{\partial f}{\partial y}(x,y_j)\eta_j = & \left(
\frac{\partial^2 L}{\partial y^2}(x,y_j)
-\varphi_j \frac{\partial^2 f}{\partial y^2}(x,y_j)\right)z_j\text{ in }\Omega,\, \eta_j = 0\text{ on }\Gamma
\end{align*}
\;
Solve the quadratic problem\label{line9}
\[
\mathrm{(}Q_j\mathrm{)}\qquad \min_{v\in L^2(\mathbb{J}_j)} H_j(v) := \frac{1}{2} F''(u_j)(\chi_{{}_{\mathbb{J}_j}} v)^2 -\int_{\mathbb{J}_j}(\kappa w_j -\eta_j) v\mathrm{d}x
\]
Name $v_{\mathbb{J}_j}$ its solution.\;
Set $v_{j} =  \chi_{{}_{\mathbb{A}_j}}w_{j} + \chi_{{}_{\mathbb{J}_j}} v_{\mathbb{J}_j}$\;
Set $u_{j+1}=u_j+v_j$ and $j=j+1$.\;
}
\end{algorithm2e}

\section{Some computational details and numerical examples}
\label{S4}
Let us comment on how to solve the quadratic problem $\mathrm{(}Q_j\mathrm{)}$ that appears in line \ref{line9} of Algorithm \ref{Alg2}. Notice that we can write $H_j(v) = \frac{1}{2}(v,A_j v)_{L^2(\mathbb{J}_j)} - (b_j,v)_{L^2(\mathbb{J}_j)}$, where $b_j = \chi_{{}_{\mathbb{J}_j}} (\kappa w_{j} -\eta_{j})$ and we can compute $A_j v$ using Algorithm \ref{Alg3}. Therefore $\mathrm{(}Q_j\mathrm{)}$ can be solved using, e.g., the conjugate gradient method without need of the explicit computation of the Hessian $F''(u_j)$.
\LinesNumbered
\begin{algorithm2e}[h!]
\caption{Computation of the product Hessian-vector}\label{Alg3}
\DontPrintSemicolon
Solve
$A z +\frac{\partial f}{\partial y}(x,y_j)z = \chi_{{}_{\mathbb{J}_j}}v \text{ in }\Omega,\ z=0\text{ in }\Gamma$\;
Solve
$A^*\eta +\frac{\partial f}{\partial y}(x,y_j)\eta = \left(
\frac{\partial^2 L}{\partial y^2}(x,y_j)
-\varphi_j \frac{\partial^2 f}{\partial y^2}(x,y_j)\right)z\text{ in }\Omega,\ \eta=0\text{ in }\Gamma$\;
Set $A_jv = \chi_{{}_{\mathbb{J}_j}}(\eta+\kappa v)$\;
\end{algorithm2e}

From the computational point of view, at each step of Algorithm \ref{Alg2} we have to solve one non-linear partial differential equation and several linear partial differential equations: three before solving the quadratic problem and two at each step of the conjugate gradient method that we use to solve the quadratic problem. When discretized, all these linear equations share either the same coefficient matrix (or its transpose in the case of a nonsymmetric problem; see e.g. \cite{CMR2020}). Therefore, advantage can be taken from a single factorization. If the nonlinear equation at iteration $j+1$ is solved using Newton's method, the matrix of the linear problem to be solved in the first iteration in this subproblem is the same as the matrix used for the linear equations at iteration $j$.

We present one example posed in a 2D domain and another one in a 3D domain.
To solve the problem we use the finite element approximation studied in \cite{CHW2012-NM}: the state, the adjoint state and the control are discretized using continuous piecewise linear elements and the Tikhonov and sparsity terms are discretized using the composite trapezoid formula. We stop the algorithm when $\delta_j = \frac{\| v_j\|_{L^2(\Omega)}}{\max\{1,\| u_{j+1}\|_{L^2(\Omega)}\}} <  5\times 10^{-14}$ or when $J(u_j)$ and $J(u_{j+1})$ are equal up to machine precision. At each iteration, the solution of the quadratic subproblem $(Q_j)$ is obtained with Matlab built-in command \texttt{pcg} and the nonlinear equation in line \ref{line3alg2} is solved using Newton's method. The tolerance  $5\times 10^{-14}$ is used for both subproblems.

\paragraph*{Example 1} We consider the data of Example 1 in \cite{CHW2012-NM}, where convergence of the finite element approximation of \Pb is studied and error estimates in terms of the discretization parameter are obtained. $\Omega = B_1(0,0)\subset \mathbb{R}^2$, $f(x,y)=y^3$, $L(x,y)=\frac{1}{2}(y-y_d(x))^2$ with $y_d(x) = 3\sin(2\pi x_1)\sin(\pi x_2)\mathrm{e}^{x_1}$, $\kappa =0.002$, $\gamma = 0.03$, $\alpha = -12$, and $\beta = 12$.

As in \cite{CHW2012-NM}, we solve the problem in a mesh of size $h= 2^{-7}$ ($1.3\times 10^5$ elements, 66049 nodes). In order to get an initial point $u_0$ close enough to $\bar u$, we take as $u_0$ the solution of the discretized problem with mesh size $h=2^{-6}$.

We have summarized the convergence history in Table \ref{T1}.
The superlinear order of convergence can be appreciated in the way the order of magnitude of the error between iterations $\delta_j$ vary in the first steps: $-2$, $-5$, $-7$, $-14$. We find numerically that $|\mathbb{J}| = 0.678$, $|\mathbb{A}^\beta|=0.310$,      $|\mathbb{A}^\alpha|=0.310$,  $|\mathbb{A}^0|=1.844$ and
$|\Sigma_{\bar u}| =0$.

\begin{table}[h!]
  \centering
  \begin{tabular}{ccccc}
    $j$ &       $J(u_j)$            &   $\delta_j$   &  $\sharp$Newton & $\sharp$CG  \\ \hline
  0&  11.141742584195615 &   $1.5\times 10^{-2}$  &   4 &  13      \\
  1&  11.141687025807151 &   $4.5\times 10^{-5}$  &   3 &  12      \\
  2&  11.141686904484867 &   $1.0\times 10^{-7}$  &   3 &  13      \\
  3&  11.141686904484866 &   $2.4\times 10^{-14}$  &   2 &  14      \\
  4 & 11.141686904484862 &                        &   1  &
  \end{tabular}
  \caption{Convergence history of the problem in the example. $\sharp$Newton is the number of Newton iterations to solve the nonlinear PDE in line \ref{line3alg2} and $\sharp$CG is the number of iterations of the conjugate gradient method used to solve $(Q_j)$ in Algorithm \ref{Alg2}}\label{T1}
\end{table}

\paragraph*{Example 2}
Consider $\Omega=(0,1)^3\subset\mathbb{R}^3$,  $f(x,y)= |y|^3 y$, $L(x,y)=\frac{1}{2}(y-y_d(x))^2$ with $y_d = \prod_{i=1}^{3} 8 x_i(1-x_i)$, $\kappa =0.1$, $\gamma = 0.05$, $\alpha = -1$, $\beta = 1$.

 We use a mesh of size $h=2^{-5}$ ($1.97\times 10^5$ elements, 35937 nodes) and start with $u_0=y_d$.

   We have summarized the convergence history in Table \ref{T2}.
   The {superlinear} order of convergence can be appreciated in the way the order of magnitude of the error between iterations $\delta_j$  vary in the first steps: $0$, $-2$, $-7$, $-16$. We find numerically that $|\mathbb{J}| = 0.323$, $|\mathbb{A}^\beta|=0.157$,      $|\mathbb{A}^\alpha|=0$,  $|\mathbb{A}^0|=0.520$ and
$|\Sigma_{\bar u}| =0$.

\begin{table}[h!]
  \centering
  \begin{tabular}{ccccc}
    $j$ &       $J(u_j)$            &   $\delta_j$   &  $\sharp$Newton & $\sharp$CG  \\ \hline
  0&  5.1160436513941248 &   $2.6\times 10^{0}$  &   4 &   4 \\
  1&  4.8088004565179974 &   $1.0\times 10^{-2}$  &   4 &   4 \\
  2&  4.8087950298698070 &   $2.7\times 10^{-7}$  &   3 &   4 \\
  3&  4.8087950298698035 &   $6.6\times 10^{-16}$  &   2 &   5 \\
  4 & 4.8087950298698035 &                         &  1  &

  \end{tabular}
  \caption{Convergence history of the problem in the example 2. $\sharp$Newton is the number of Newton iterations to solve the nonlinear PDE in line \ref{line3alg2} and $\sharp$CG is the number of iterations of the conjugate gradient method used to solve $(Q_j)$ in Algorithm \ref{Alg2}}\label{T2}
\end{table}


\end{document}